\documentclass[12pt]{article}

\usepackage{amsgen,amsmath,amstext,amsbsy,amsopn,amsfonts,amssymb}
\usepackage[dvips]{graphicx}

\title{Some properties of index of Lie algebras}
\author{Vladimir Dergachev}

\def\buzz#1{{\normalfont{\bf{[#1]}}}}

\def\ind{\textrm{\normalfont ind\,}}

\def\acr{\nonumber \\}

\def\rank{\textrm{\normalfont rank\,}}

\def\stab{\textrm{\normalfont Stab\,}}
\def\ad#1{\textrm{\normalfont ad}_{{#1}}\,}
\def\Ad#1{\textrm{\normalfont Ad}_{{#1}}\,}
\def\coad#1{\textrm{\normalfont ad}^*_{{#1}}\,}
\def\coAd#1{\textrm{\normalfont Ad}^*_{{#1}}\,}
\def\diag{\textrm{\normalfont diag}}
\def\Mat{\textrm{\normalfont Mat}}

\def\Stab{\textrm{\normalfont Stab}}
\def\stab{\Stab}

\def\gl{{\mathfrak{gl}}}
\def\imply{\Rightarrow}
\newtheorem{theorem}{Theorem}
\newtheorem{definition}{Definition}
\newtheorem{example}{Example}
\newtheorem{proposition}[theorem]{Proposition}
\newtheorem{lemma}[theorem]{Lemma}

\newenvironment{proof}{\begin{trivlist}\item[\hskip%
\labelsep{\bf Proof\quad}]}%
{\hfill\qed\end{trivlist}}
\newcommand{\qed}{{\unskip\nobreak\hfil\penalty50\hskip .001pt \hbox{}
          \nobreak\hfil
          \vrule height 1.2ex width 1.1ex depth -.1ex
           \parfillskip=0pt\finalhyphendemerits=0\medbreak}\rm}

\makeatletter
\newcount\@refno
\def\nextref#1#2#3#4{\advance\@refno\@ne
\if@filesw \immediate\write\@auxout
   {\string\bibcite{#1}{\number\@refno}}\fi   
    {{\number\@refno}.\quad{ #2}, {#3}, { #4}.\hfill\\}}

\newcommand{\references}{
\section*{References}
\frenchspacing
\entries\par}

\newcommand{\entries}{
\nextref{DK}
{Vladimir Dergachev and Alexandre Kirillov}
{\em Index of Lie algebras of Seaweed type}
{to appear in {\em Journal of Lie theory}}
\nextref
{D}
{J.Dixmier}
{``Enveloping algebras''}
{American Mathematical Society, Providence (1996) 1-379 }
\nextref{E1}{Elashvili, A. G}
{\em Frobenius Lie algebras}
{Funktsional. Anal. i Prilozhen. {\bf 16} (1982), 94--95}
\nextref{E2}{---}
{\em On the index of a parabolic subalgebra 
of a semisimple Lie algebra}
{Preprint, 1990}
\nextref{GK}
{I.M.Gelfand and A.A.Kirillov}
{\em Sur les corps li\'es aux alg\'ebres enveloppantes des alg\'ebres de Lie}
{ Publications math\'ematiques {\bf 31} (1966) 5-20}

}

\begin{document}
\maketitle
\tableofcontents

\section{Introduction}
In paper \cite{DK} the author and Alexandre Kirillov have computed index
\footnote{For the definition of index of Lie algebra see \cite{D}. }
for a family of subalgebras of $\gl(n)$. The papers \cite{E1}, \cite{E2} and some
computations done by the author provide insight into situation with subalgebras
of other classical groups. One interesting property of this numeric data is that
while subalgebras of $\gl(n)$ exhibit much variety in the possible values of 
the index the subalgebras of other simple groups generally do not. We believe
that this is due to the fact that Lie bracket on $\gl(n)$ can be derived from
multiplication in the associative algebra of matrices, while groups from other
series do not possess this property. Thus one expects to discover that index
will exhibit special characteristics in relation to Lie algebras derived from associative
algebras. In this paper we explore this idea. We discover that index possesses
certain "convexity" properties with respect to the operation of tensor product
of associative algebras. Moreover there is a large family of associative algebras
for which the convexity inequalities become precise thus shedding light on 
the richness of structure observed among subalgebras of $\gl(n)$.

\section{Index of Lie algebras}
\begin{definition} \buzz{Lie algebra}

A Lie algebra $\mathfrak g$ over field $k$ is a vector space over $k$ with operation
$[\cdot,\cdot]:\mathfrak g \times \mathfrak g \rightarrow \mathfrak g$ that satisfies
the following properties:
\begin{enumerate}
\item skew-symmetry:
$$
[a,b]=-[b,a]
$$
\item Jacobi identity:
$$
[[a,b],c]+[[c,a],b]+[[b,c],a]=0
$$
\end{enumerate}
\end{definition}

\begin{definition} \buzz{Index of Lie algebra}
Pick a basis $e_i\in \mathfrak g$. Let $B$ be the multiplication matrix of the
bracket product: $b_{i,j}=[e_i,e_j]$. $B$ is then skew-symmetric and it's elements
can be considered as polynomials over $\mathfrak g^*$ - dual space to $\mathfrak g$.

For $f\in \mathfrak g^*$ we define $\left.B\right|_f$ to be the result of evaluation
of elements of $B$ on $f$, that is $\left.B\right|_f=\left\|f([e_i,e_j])\right\|$.

By definition the index of Lie algebra $\mathfrak g$ is
$$
\ind \mathfrak g=\min_{f\in\mathfrak g^*} \dim \ker \left.B\right|_f
$$
\end{definition}
 
We will use $\Omega$ to denote an element of $S(\mathfrak g^*)\otimes_k \wedge^2\mathfrak g$
that corresponds to $B$.
Let $r$ be the maximum number such that $\wedge^r\Omega\neq 0$. 

\begin{proposition}
\label{BiOmegaProposition}
The following numbers are equal:
\begin{itemize}
\item $\ind \mathfrak g$
\item $\dim \mathfrak g - 2r$
\item $\dim \ker B$ , where $B$ is considered as a matrix with coefficients
in the field $k(\mathfrak g^*)$ of rational functions on $\mathfrak g^*$
\item $\dim \ker B_f$ for generic (in Zariski sense) $f\in \mathfrak g^*$ 
\end{itemize}
\end{proposition}

For more extended discussion of $B$ and $\Omega$
see \cite{DK}.

\section{Tensor products}
\begin{definition} \buzz{Tensor product of matrices}

Let $A$ and $B$ be two matrices
with coefficients in rings $\mathcal R_1$ and $\mathcal R_2$
respectively. Let ring $\mathcal R$ have the property that $\mathcal R \subset \mathcal R_1$
and $\mathcal R \subset \mathcal R_2$.
The tensor product $A \otimes_{\mathcal R} B$ is defined as a block matrix with each block $(i,j)$ having dimensions of matrix $B$
and equal to $A_{i,j}\otimes_{\mathcal R} B$, that is the matrix obtained from $B$ by taking tensor products of
a certain element of $A$ with entrees of $B$. Thus $A \otimes B$ has
coefficients in $\mathcal R_1 \otimes_{\mathcal R} \mathcal R_2$.

\end{definition}

\begin{proposition}The tensor product of matrices has the following
properties:
\begin{enumerate}
\item distributive w.r.t. addition
\item $\left(A\otimes_{\mathcal R}B\right)\cdot \left(C\otimes_{\mathcal R}D\right)=
\left(AC\right)\otimes_{\mathcal R}\left(BD\right)$, where ($A$,$C$) and ($B$,$D$) are
pairs of matrices of the same shape
\item $\left(A\otimes_{\mathcal R}B\right)^{-1}=\left(A^{-1}\right)\otimes_{\mathcal R}\left(B^{-1}\right)$
\end{enumerate}
\end{proposition}

\begin{theorem}
Let $A$ and $B$ be square matrices of dimensions $k$ and $n$ respectively,
with coefficients in commutative rings $\mathcal R_1$ and $\mathcal R_2$. Let
ring $\mathcal R$ have the property that $\mathcal R \subset \mathcal R_1$
and $\mathcal R \subset \mathcal R_2$. Then
$$
\det\left(A\otimes_{\mathcal R}B\right)=\left(\det
A\right)^n\otimes_{\mathcal R}\left(\det
B\right)^k
$$
\end{theorem}
\begin{proof}
\begin{enumerate}
\item If $A$ and $B$ are diagonal the statement is proved by a simple
computation.
\item Let $R_1=R_2=R=\mathbb{C}$. Let $A=C_1D_1C_1^{-1}$ and
$B=C_2D_2C_2^{-1}$ where $D_1$ and $D_2$ are diagonal. Then 
\begin{eqnarray}
\lefteqn{A\otimes_{\mathcal R}B=\left(C_1D_1C_1^{-1}\right)\otimes_{\mathcal
R}\left(C_2D_2C_2^{-1}\right)= }\nonumber\\
& & \qquad
=\left(C_1\otimes_{\mathcal R}C_2\right)
\left(D_1\otimes_{\mathcal R}D_2\right)
\left(C_1^{-1}\otimes_{\mathcal R}C_2^{-1}\right)=
 \nonumber \\
& & \qquad
=\left(C_1\otimes_{\mathcal R}C_2\right)
\left(D_1\otimes_{\mathcal R}D_2\right)
\left(C_1\otimes_{\mathcal R}C_2\right)^{-1}
 \nonumber 
\end{eqnarray}
and
\begin{eqnarray}
\lefteqn{\det\left(A\otimes_{\mathcal R}B\right)=} \nonumber\\
& &\qquad
=\det\left(\left(C_1\otimes_{\mathcal R}C_2\right)
\left(D_1\otimes_{\mathcal R}D_2\right)
\left(C_1\otimes_{\mathcal R}C_2\right)^{-1} \right)=
\nonumber \\
& &\qquad
=\det\left(
D_1\otimes_{\mathcal R}D_2
 \right)=
\nonumber \\
& & \qquad
=\left(\det D_1\right)^n \left(\det D_2\right)^k=\nonumber\\
& & \qquad
=\left(\det A\right)^n \left(\det B\right)^k\nonumber
\end{eqnarray}
\item Since both sides of the equation $\det\left(A\otimes_{\mathcal R}B\right)=\left(\det
A\right)^n\otimes_{\mathcal R}\left(\det
B\right)^k$ are polynomials in elements of $A$ and $B$ with integral coefficients and we know that over $\mathbb{C}$
all generic $A$ and $B$ satisfy the equation we must have that the polynomials
are identical. This concludes the proof of the theorem.
\end{enumerate}
\end{proof}

\begin{theorem}
\label{MatTensProduct}Let $\mathcal F$, $\mathcal F_1$ and $\mathcal F_2$ be
three fields, such that $\mathcal F\subset \mathcal F_1$ and $\mathcal
F\subset \mathcal F_2$. Let $A_1$ and $B_1$ be two matrices with
coefficients in $\mathcal F_1$ and $A_2$ and $B_2$ be two matrices with
coefficients in $\mathcal F_2$.

Then
\begin{enumerate}
\item $\dim \ker \left({B_1 \otimes_{\mathcal F} A_2 + A_1 \otimes_{\mathcal F} B_2}\right)\ge 
\dim \ker B_1 \cdot \dim \ker B_2$
\item For generic $A_1$ and $A_2$ (in Zariski topology) the inequality above is precise.
\end{enumerate}
\end{theorem}
\begin{proof}

\begin{enumerate}
\item Using elementary linear algebra one obtains matrices $C_1$ and
$C_2$, such that $\rank C_1=\dim \ker B_1$ and $\rank C_2=\dim
\ker B_2$ and also $B_1C_1=0$ and $B_2C_2=0$.

Because of properties of tensor product we have 
$$
\rank \left(C_1 \otimes_{\mathcal F} C_2\right)=\rank C_1 \cdot \rank C_2
$$

And since 
\begin{eqnarray}
\lefteqn{ \left({B_1 \otimes_{\mathcal F} A_2 + A_1 \otimes_{\mathcal F}
B_2}\right) \cdot \left(C_1 \otimes_{\mathcal F} C_2\right)=}\acr
& &\qquad 
=\left(B_1C_1\right) \otimes_{\mathcal F} \left(A_2C_2\right)+
\left(A_1C_1\right)\otimes_{\mathcal F}\left(B_2C_2\right)=0\nonumber
\end{eqnarray}
we have $\dim \ker \left({B_1 \otimes_{\mathcal F} A_2 + A_1 \otimes_{\mathcal F} B_2}\right)\ge 
\rank \left(C_1\otimes_{\mathcal F}C_2\right)$ which concludes the first
part of the proof.

\item Let $\hat{A_1}=\left\|a^1_{i,j}\right\|$ and
$\hat{A_2}=\left\|a^2_{i,j}\right\|$ where $\left\{a^1_{i,j}\right\}$ and
$\left\{a^2_{i,j}\right\}$ are two families of independent variables.
Let $\hat{\mathcal F_1}={\mathcal F_1}\left(\left\{a^1_{i,j}\right\}\right)$
and $\hat{\mathcal F_2}={\mathcal F_2}\left(\left\{a^2_{i,j}\right\}\right)$
be fields of rational functions over $\left\{a^1_{i,j}\right\}$ and
$\left\{a^2_{i,j}\right\}$ correspondingly.

Let 
$$
\vec{W}=v_1 \otimes_{\mathcal F} w_1+...+v_n \otimes_{\mathcal F}
w_n
$$
 be a vector from $\ker \left(B_1 \otimes_{\mathcal F}\hat{A_2} + \hat{A_1}
\otimes_{\mathcal F}B_2\right)$ ($v_i$ have coefficients in $\hat{\mathcal F_1}$ and
$w_j$ have coefficients in $\hat{\mathcal F_2}$).
 Suppose that decomposition of $\vec{W}$
above is simple in the sense that neither two of $v_i$ ($w_j$) are
proportional and $n$ is the minimum possible number for such a
decomposition. Since $\vec{W}$ is defined over $\hat{\mathcal
F_1}\otimes_{\mathcal F}\hat{\mathcal F_2}$ we can find a polynomial $p$
from 
$$
{\mathcal Q}={\mathcal
F_1\left[\left\{a^1_{i,j}\right\}\right]}\otimes_{\mathcal F}{\mathcal
F_2\left[\left\{a^2_{i,j}\right\}\right]}
$$
such that all elements of $p\vec{W}$ are from $\mathcal Q$. We now
introduce a bi-grading in $\mathcal Q$: 
\begin{eqnarray}
\deg {\mathcal F}&=&(0,0) \nonumber \\
\deg {a^1_{i,j}}&=&(1,0)  \nonumber \\
\deg {a^2_{i,j}}&=&(0,1)  \nonumber 
\end{eqnarray}

The matrix $B_1 \otimes_{\mathcal F}\hat{A_2} + \hat{A_1}
\otimes_{\mathcal F}B_2$ has thus two parts - of degree $(0,1)$ and
$(1,0)$. Therefore, every element of it's kernel is the tensor product 
of an element of $\ker B_1$ and an element of $\ker B_2$ - with possible
coefficient from $\hat{\mathcal F_1}\otimes_{\mathcal F}\hat{\mathcal
F_2}$. This proves the statement for the particular case of $\hat{A_1}$
and $\hat{A_2}$, thus implying that the equality holds for almost all
$A_1$ and $A_2$ (in Zariski sense).
\end{enumerate}
\end{proof}

\section{Associative algebras}
\begin{definition}

Let $\mathfrak A$ be an associative algebra. Then $\mathfrak A^L$
denotes a Lie algebra with the bracket defined as $[a,b]=ab-ba$.
\end{definition}

\begin{theorem}\buzz{Convexity property of index}
Let $\mathfrak A$ and $\mathfrak B$ be two finite dimensional associative 
algebras over field $\mathcal F$. Then
\begin{equation}
\ind \left({\mathfrak A \otimes_{\mathcal F} \mathfrak B}\right)^L\ge \ind {\mathfrak A}^L \cdot \ind
{\mathfrak B}^L
\label{tensorinequality}
\end{equation}
\end{theorem}
\begin{proof}
Let $\left\{e_i\right\}$ be a basis in $\mathfrak A$ and $\left\{g_j\right\}$
be a basis in $\mathfrak B$. Denote by $A_1$ the multiplication matrix of $\mathfrak A$,
that is the matrix with coefficients in (first degree) polynomials over $\mathfrak A^*$,
element $(i,j)$ of $A_1$ is equal to $e_i \cdot e_j$. Let $A_2$ denote the multiplication
matrix of $\mathfrak B$.

The multiplication matrix $B_1$ of {\em Lie algebra} $\mathfrak A^L$ is given
by the formula $B_1=A_1-A_1^T$  (here $(\cdot)^T$ denotes transposition).
$B_2=A_2-A_2^T$ gives the multiplication matrix of $\mathfrak B^L$.

The multiplication matrix of $\mathfrak A \otimes_{\mathcal F} \mathfrak B$ in the 
basis $\left\{e_i\otimes_{\mathcal F}g_j\right\}$ is $A_1 \otimes_{\mathcal F}A_2$.
The multiplication matrix $B$ for $\left(\mathfrak A \otimes_{\mathcal F} \mathfrak B\right)^L$
is computed as follows:
\begin{eqnarray}
\lefteqn{B=A_1 \otimes_{\mathcal F}A_2-\left(A_1 \otimes_{\mathcal F}A_2\right)^T=}\nonumber \\
&&\quad=A_1 \otimes_{\mathcal F}A_2-A_1^T \otimes_{\mathcal F}A_2^T=   \nonumber \\
&&\quad=A_1 \otimes_{\mathcal F}A_2-A_1^T \otimes_{\mathcal F}A_2+
A_1^T \otimes_{\mathcal F}A_2-A_1^T \otimes_{\mathcal F}A_2^T=   \nonumber \\
&&\quad=B_1 \otimes_{\mathcal F}A_2+A_1^T \otimes_{\mathcal F}B_2    \nonumber
\end{eqnarray}
Theorem \ref{MatTensProduct} implies that $\dim \ker B\ge \dim \ker B_1 \cdot
\dim \ker B_2$. 
One concludes the proof by applying proposition \ref{BiOmegaProposition}.
\end{proof}

\begin{example}\normalfont
Let $\mathfrak A=a\mathbb C+b\mathbb C$ be a 2-dimensional algebra over $\mathbb C$ (you can replace $\mathbb C$ with
you favorite field) with the following
multiplication table:
\begin{center}
\begin{tabular}{c|cc}
$\times_{\mathfrak A}$& $a$ & $b$\\
\hline
$a$ & $a$ & $b$ \\
$b$ & $a$ & $b$ \\
\end{tabular}
\end{center}
Let $\mathfrak B$ be an arbitrary associative algebra over $\mathbb C$ with multiplication 
table $A$. The the multiplication table of $\mathfrak A\otimes_{\mathbb C} \mathfrak B$ is:
\begin{center}
\begin{tabular}{c|cc}
$\times_{\mathfrak A\otimes_{\mathbb C} \mathfrak B}$& $a\otimes_{\mathbb C} \mathfrak B $& $b\otimes_{\mathbb C} \mathfrak B$\\
\hline
$a\otimes_{\mathbb C} \mathfrak B$ & $aA$ & $bA$ \\
$b\otimes_{\mathbb C} \mathfrak B$ & $aA$ & $bA$ \\
\end{tabular}
\end{center}
The index of algebra $\mathfrak A$ is equal to $0$. The index of algebra $\mathfrak A\otimes_{\mathbb C}\mathfrak B$
can be computed as follows:
\item
$$
B_{\mathfrak A\otimes_{\mathbb C}\mathfrak B}=\left(
\begin{matrix}
aA-aA^T & bA-aA^T \\
aA-bA^T & bA-bA^T
\end{matrix}
\right)
$$
By definition $\ind \mathfrak A\otimes_{\mathbb C}\mathfrak B=\ker 
B_{\mathfrak A\otimes_{\mathbb C}\mathfrak B}$. Thus we want to find all 
solutions $(v_1,v_2)$ (in rational functions on $\left(\mathfrak A\otimes_{\mathbb C}\mathfrak B\right)^*$)
of the following equations:
$$
\left\{
\begin{matrix}
\left(aA-aA^T\right)v_1 +\left(bA-aA^T\right)v_2&=&0 \nonumber \\
\left(aA-bA^T\right)v_1 +\left(bA-bA^T\right)v_2&=&0  \nonumber
\end{matrix}
\right.
$$
A few straightforward transformations lead us to the following system:
$$
\left\{
\begin{matrix}
(a-b)A^T\left(v_1+v_2\right)&=&0 \nonumber \\
(a-b)A\left(av_1+bv_2\right)&=&0  \nonumber
\end{matrix}
\right.
$$
Since $a$ and $b$ are invertible and independent the index of $\mathfrak A\otimes_{\mathbb C}\mathfrak B$
is equal to twice the dimension of the kernel of $A$. 
\qed
\end{example}

In the particular case considered above the inequality \ref{tensorinequality} is
precise only when the algebra $\mathfrak B$ has non-degenerate multiplication table.

\section{Characteristic polynomial of associative algebra}
\begin{definition}\buzz{Multiplicative group}
Let $\mathfrak A$ be an associative algebra with unity. Then $\mathfrak A^*$ -
the set of all invertible elements of $\mathfrak A$ - can be regarded as a group
under the operation of multiplication.

The Lie algebra of $\mathfrak A^*$ is the algebra $\mathfrak A$ itself with
commutator 
$$
[a,b]=ab-ba
$$
\end{definition}

\begin{definition}\buzz{Adjoint action}
Let $g\in\mathfrak A^*$ and $Y\in\mathfrak A={\left(\mathfrak
A^*\right)}^L$.  By definition {\it adjoint action} of $g$ on $X$ is
$$
\Ad{g} Y=gYg^{-1}
$$
Let $X\in\mathfrak A$. The adjoint action of $X$ on $Y$ is defined as
$$
\ad{X}Y=XY-YX
$$
\end{definition}

\begin{definition}\buzz{Coadjoint action}
Let $g\in\mathfrak A^*$ and $F\in\mathfrak A'$ (here $\mathfrak A'$ is
the space of linear functionals on $\mathfrak A$). Let $Y\in\mathfrak
A$.
By definition {\it coadjoint action} of $g$ on $F$ is
$$
\left(\coAd{g} F\right)(Y)=F\left(g^{-1}Yg\right)
$$
Let $X\in\mathfrak A$. The coadjoint action of $X$ on $F$ is defined as
$$
\left(\coad{X}F\right)(Y)=F\left(-(XY-YX)\right)
$$
\end{definition}

\begin{definition}\buzz{Characteristic polynomial of associative algebra}
Let $\left\{e_i\right\}$ be some basis of $\mathfrak A$. Let $A$ be the
multiplication table of $\mathfrak A$ considered as matrix which
coefficients are polynomials over $\mathfrak A'$. Then
$$
\chi\left(\lambda,\mu,F\right)=\det\left(\lambda A+\mu A^T\right)
$$
is the {\it characteristic polynomial} of algebra $\mathfrak A$. 
Here $F$ is an element of $\mathfrak A'$.
\end{definition}

\begin{lemma}
Characteristic polynomial is defined up to a scalar multiple.
\end{lemma}
\begin{proof}
Indeed, change of basis $\left\{e_i\right\}$ replaces $A$ with $CAC^T$.
And 
$$
\det\left(\lambda CAC^T+\mu CA^TC^T\right)=\left(\det
C\right)^2\det\left(\lambda A+\mu A^T\right)
$$
\end{proof}

\begin{theorem}\label{chi_invariance}
The characteristic polynomial is quasi-invariant under coadjoint action.
That is 
$$
\chi\left(\lambda,\mu,\coAd{g}F\right)=\left(\det \Ad{g}\right)^{-2}\chi\left(\lambda,\mu,F\right)
$$
\end{theorem}
\begin{proof}
Indeed, the matrix element $(i,j)$ of $A$ evaluated in point $F$ is
given by the expression $F(e_ie_j)$. Since
$$
\left(\coAd{g}F\right)(e_ie_j)=F(g^{-1}e_ie_jg)=F((g^{-1}e_ig)(g^{-1}e_jg))
$$
the substitution $F\rightarrow \coAd{g}F$ is  equivalent to change of
basis induced by the matrix $\textrm{\normalfont{Ad}}_{g}^{-1}$. Applying the result of
previous lemma one easily obtains the statement of the theorem.
\end{proof}

\begin{theorem}\label{ext_cayley}\buzz{Extended Cayley theorem}
Let $A$ and $B$ be two $n\times n$ matrices over an algebraically closed
field $k$ and $C$ and $D$ be two $m\times m$
matrices over the same field $k$. Define $\chi(\lambda,\mu)=\det(\lambda A+\mu B)$. Then
$$
\det(\lambda A\otimes C +\mu B \otimes D)=det(\chi(\lambda C,\mu D))
$$
\end{theorem}

Before proceeding with the proof we must explain in what sense we 
consider $\chi(\lambda C,\mu D)$. Indeed, matrices $C$ and $D$ might
not commute making $\chi(\lambda C,\mu D)$ ambiguous.
In our situation the right definition is as follows:

First, we notice that $\chi(\lambda,\mu)$ is homogeneous, thus it
can be decomposed into a product of linear forms ($k$ is algebraicly closed):
$$
\chi(\lambda,\mu)=\prod_i(\lambda \alpha_i +\mu \beta_j)
$$

Then we define 
$$
\chi(\lambda C,\mu D)=\prod_i(\lambda \alpha_i C +\mu \beta_j D)
$$

There is still some ambiguity about the order in which we multiply 
linear combinations of $C$ and $D$ but it does not affect the value
of $\det(\chi(\lambda C,\mu D))$.

\begin{proof}

{\bf Step 1.} Let $A$ and $C$ be identity matrices of sizes $n\times n$ and
$m \times m$ respectively. 

Then 
$$
\chi(\lambda,\mu)=\prod_i(\lambda+\mu \gamma_i)
$$
where $\gamma_i$ are eigenvalues of $B$.
\begin{eqnarray}
\det(\chi(\lambda,\mu D))=\det\left(\prod_i(\lambda+\mu \gamma_i D) \right)\qquad\cr
\qquad=\prod_i\det\left(\lambda+\mu \gamma_i D \right)=
\prod_{i,j}\left(\lambda+\mu \gamma_i \epsilon_j\right)
\end{eqnarray}
where $\epsilon_j$ are eigenvalues of $D$.

On the other hand one easily derives that eigenvalues of $B\otimes D$ are 
$\gamma_i \epsilon_j$ and thus
$$
\det(\lambda +\mu B\otimes D)=\prod_{i,j}\left(\lambda+\mu \gamma_i\epsilon_j\right)
=\det(\chi(\lambda,\mu D))
$$

{\bf Step 2.}
Let us assume now only that matrices $A$ and $C$ are invertible.

We have 
$$
\chi(\lambda,\mu)=\det(\lambda A+\mu B)=\det(A)\det(\lambda+\mu A^{-1}B)
$$

Denote $\chi'(\lambda,\mu)=\det(\lambda+\mu A^{-1}B)=\prod_i(\lambda+\mu \gamma_i)$, 
where $\gamma_i$ are eigenvalues of $A^{-1}B$.
\begin{eqnarray}
\lefteqn{\det(\chi(\lambda C,\mu D))=\det\left(\det(A)\prod_i(\lambda C+\mu \gamma_i D)\right)=} \cr
&&=\det(A)^m\prod_i\det(\lambda C+\mu \gamma_i D)=\cr
&&=\det(A)^m\prod_i\left( \det(C)\det(\lambda+\mu \gamma_i C^{-1}D)\right)=\cr
&&=\det(A)^m \det(C)^n \det(\chi'(\lambda,\mu C^{-1}D))\nonumber
\end{eqnarray}

By {\em step 1} we have 
$$
\det(\lambda+\mu (A^{-1}B)\otimes(C^{-1}D))=\det(\chi'(\lambda,\mu C^{-1}D))
$$

Observing also that $\det(A\otimes C)=\det(A)^m\det(C)^n$ we obtain
\begin{eqnarray}
\det(\chi(\lambda C,\mu D))=\det(A\otimes C)\det(\lambda+\mu (A^{-1}B)\otimes(C^{-1}D))=\cr
\qquad =\det(\lambda A\otimes C+\mu B\otimes D)\nonumber
\end{eqnarray}
which is the desired formula.

{\bf Step 3.} To prove the formula in general we observe that the both parts
involve only polynomials in entries of matrices $A$,$B$,$C$ and $D$. Since
the restriction that $A$ and $C$ be invertible selects a Zariski open subset,
the formula should hold for all $A$,$B$,$C$,$D$ by continuity.
\end{proof}

\begin{proposition}
\label{p_divisibility}
Let $A$ and $B$ be two $n\times n$ matrices. Let $K=\dim \ker A$ and $M=\dim \ker B$.
 Then $\det(\lambda A+\mu B)$ is 
divisible by $\lambda^M\mu^K$.
\end{proposition}

\begin{theorem}\label{upper_bound_theorem}
Let $\mathfrak A$ be an associative algebra. Then it's characteristic polynomial
is divisible by 
$$
(\lambda\mu)^{\dim \ker A}(\lambda+\mu)^{\ind \mathfrak A}
$$
\end{theorem}
\begin{proof}
Indeed, using proposition \ref{p_divisibility} and the fact that
$\dim \ker A=\dim \ker A^T$ the characteristic polynomial
$\det(\lambda A + \mu A^T)$ is divisible by $(\lambda \mu)^{\dim \ker A}$.

And in view of 
$$
\det(\lambda A + \mu A^T)=\det((\lambda+\mu)A-\mu(A-A^T))
$$
and
$$
\dim \ker (A-A^T)=\ind \mathfrak A
$$
it is divisible by $(\lambda+\mu)^{\ind \mathfrak A}$.
\end{proof}

\begin{example}\normalfont \buzz{$\textrm{GL}(2)$}
Let $a$,$b$,$c$,$d$ denote the matrix units $E_{1,1}$,$E_{1,2}$,$E_{2,1}$ and
$E_{2,2}$ correspondingly. Then the multiplication table $A$ is 
$$
\begin{array}{cccccc}
&\vline & a & b & c & d \\
 \hline
a& \vline& a & b & 0  & 0 \\
b& \vline& 0 & 0 & a & b  \\
c&\vline& c  & d  & 0 & 0 \\
d& \vline& 0 & 0  & c   & d \\
\end{array}
$$

The characteristic polynomial is equal to
\begin{eqnarray}
\lefteqn{
\chi(\lambda,\mu)=\det(\lambda A+\mu A^T)=}\cr
& &=-(\lambda+\mu)^2(ad-bc)(
(\lambda-\mu)^2(ad-bc)+\lambda\mu(a+d)^2))
\end{eqnarray}
We see the absence of factors $\lambda\mu$ which corresponds to the fact that
multiplication table $A$ is non-degenerate ($\det A=-(ad-bc)^2$) and the degree
of the factor
$\lambda+\mu$ is exactly equal to the index of $\textrm{GL}(2)$.
\end{example}

\begin{example}\normalfont\label{seaweed21x12} Let $\mathfrak A$ be the following 
subalgebra of $\Mat_3(\mathbb C)$:
$$
\left(
\begin{matrix}
a & b & 0 \\
0 & c & 0 \\
0 & d & e
\end{matrix}
\right)
$$
The multiplication table $A$ of $\mathfrak A$ is 
$$
\begin{array}{ccccccc}
&\vline & a & b & c & d & e \\
\hline 
a & \vline & a & b & 0 & 0 & 0\\
b & \vline & 0 & 0 & b & 0 & 0\\
c & \vline & 0 & 0 & c & 0 & 0\\
d & \vline & 0 & 0 & d & 0 & 0\\
e & \vline & 0 & 0 & 0 & d & e
\end{array}
$$
The characteristic polynomial is equal to
$$
\chi(\lambda,\mu)=\lambda^2\mu^2(\lambda+\mu)b^2d^2(a+c+e)
$$
The degree of the factor $\lambda+\mu$ is equal to the index of $\mathfrak A$.
The presence of factor $\lambda^2\mu^2$ corresponds to the fact that multiplication
table is degenerate.
\end{example}

Algebra $\mathfrak A$ from the example \ref{seaweed21x12} facilitates construction
of the algebra with null characteristic polynomial. Indeed, by theorem \ref{ext_cayley}
characteristic polynomial of algebra $\mathfrak A\otimes \mathfrak A$ is equal to
$$
\chi_{\mathfrak A \otimes \mathfrak A}(\lambda,\mu)
=\det((\lambda A)^2(\mu A^T)^2(\lambda A+\mu A^T)b^2d^2(a+c+e))=0
$$

\section{Characteristic polynomial of $\Mat_n$}
\begin{definition}\buzz{Generalized resultant} Let $p(x)$ and $q(x)$ be two
polynomials over an algebraicly closed field. We define 
{\em generalized resultant} of $p(x)$ and $q(x)$ to be
$$
R(\lambda,\mu)=\prod_{i,j}{(\lambda \alpha_i+\mu \beta_j)}
$$
where $\left\{\alpha_i\right\}$ and $\left\{\beta_j\right\}$ are roots of
polynomials $p(x)$ and $q(x)$ respectively.
\end{definition}

Generalized resultant is polynomial in two variables. It is easy to show that
its coefficients are polynomials in coefficients of $p(x)$ and $q(x)$ so the
condition on base field to be algebraicly closed can be omitted.

\begin{theorem}The characteristic polynomial $\chi(\lambda,\mu,F)$ for algebra
$\Mat_n$ in point $F\in \Mat_n^*$ is equal to the generalized resultant of 
characteristic polynomial of $F$ (as a matrix) with itself times $(-1)^{\frac{n(n-1)}{2}}$. That is
$$
\chi(\lambda,\mu,F)=(-1)^{\frac{n(n-1)}{2}}\det(F)(\lambda+\mu)^n\prod_{i\neq j}{(\lambda \alpha_i+\mu \alpha_j)}
$$
where $\alpha_i$ are eigenvalues of $F$ (this formulation assumes
that  the base field is algebraicly closed).
\end{theorem}
\begin{proof}
We will make use of theorem \ref{chi_invariance}. The coadjoint
action on $\Mat_n^*$ is simply conjugation by invertible matrices. The generic
orbit consists of diagonalizable matrices. Thus we can compute $\chi(\lambda,\mu)$
by assuming first that $F$ is diagonal and then extrapolating the resulting polynomial
to the case of all $F$.

Assume the base field to be $\mathbb C$.
Let $F=\diag(\alpha_1,...,\alpha_n)$. We choose a basis $\left\{E_{i.j}\right\}$
of matrix units in the algebra $\Mat_n$. The only case when $F(E_{i,j}E_{k,l})$
is non-zero is when $i=l$ and $j=k$. Thus the multiplication table $A$ (restricted
to subspace of diagonal matrices in $\Mat_n^*$) is
$$
\begin{array}{ccccccccccccccccc}
& \vline &  & E_{i,i}   &  &\vline& & E_{i,j}^+ & &\vline & & E_{i,j}^- & \cr
\hline
& \vline& \alpha_1  & & 0& \vline &&&& \vline \cr
E_{i,i} & \vline&  & \ddots  & &\vline& & 0 & &\vline & & 0\cr
& \vline& 0&  & \alpha_n &\vline& & & & \vline & \cr
\hline
& \vline& & && \vline & & & & \vline &\alpha_{j'} & & 0 \cr
E_{i,j}^+ & \vline & & 0 & &\vline & & 0 & &\vline &&\ddots	\cr
& \vline& & && \vline & & & & \vline &0 & &\alpha_{j''}\cr
\hline
& \vline& & && \vline &\alpha_{i'} & &0 & \vline & \cr
E_{i,j}^- &\vline & & 0 & &\vline & &\ddots & &\vline & & 0\cr
& \vline& & && \vline &0 & &\alpha_{i''} & \vline & \cr
\end{array}
$$
here $E_{i,j}^+$ denotes elements $E_{i,j}$ with $i>j$ and $E_{i,j}^-$ denotes
elements $E_{i,j}$ with $i<j$. The matrix $\lambda A+\mu A^T$ will have 
$(\lambda+\mu)\alpha_i$ in the $E_{i,i}\times E_{i,i}$ block, and the pair 
$(E_{i,j}^+,E_{i,j}^-)$ will produce a $2\times 2$ matrix
$$
\left(
\begin{array}{cc}
0 & \lambda\alpha_j+\mu\alpha_i \cr
\lambda\alpha_i+\mu\alpha_j & 0 \cr
\end{array}
\right)
$$
Computing the determinant yields
$$
(-1)^{\frac{n(n-1)}{2}}(\lambda+\mu)^n\prod_i\alpha_i\prod_{i\neq j}{(\lambda\alpha_i+\mu\alpha_j)}
=(-1)^{\frac{n(n-1)}{2}}\prod_{i,j}{(\lambda\alpha_i+\mu\alpha_j)}
$$
thus proving the theorem for the case when $F$ is diagonal. But characteristic
polynomial $\det(F-x)$ is invariant under coadjoint action. Thus the expression
is true for all $F$ up to a possibly missing factor depending only on $F$ (but not
$\lambda$ or $\mu$) which is quasi-invariant under co\-adjoint action. However,
in view of the fact that
this multiple must be polynomial in $F$ and that the degree of the 
expression above in $F$ is exactly $n^2$ this multiple must be trivial.

The case of arbitrary field is proved by observing that both sides of the equality
are polynomials with integral coefficients and thus if equality holds over $\mathbb C$
it should hold over any field.
\end{proof}

We observe that the maximal degree of factor $(\lambda+\mu)$ in characteristic
polynomial of $\Mat_n$ is equal exactly to the index of $\Mat_n$. Moreover, in the
factorization of $\chi(\lambda,\mu)$ all factors except $(\lambda+\mu)$ depend
upon $F$ non-trivially. This allows us to prove the following theorem:
\begin{theorem}
Let $\mathfrak A$ be a an algebra with the property that the maximal $n$ for 
which $(\lambda+\mu)^n$ divides $\chi_{\mathfrak A}$ is equal to $\ind \mathfrak A$.
Then 
$$
\ind\left( \Mat_N\otimes \mathfrak A\right)=N\ind \mathfrak A
$$
\end{theorem}
\begin{proof}
Indeed, from the theorem \ref{tensorinequality} we know that 
$\ind\left( \Mat_N\otimes \mathfrak A\right)$ is at least $N\ind \mathfrak A$.
On the other hand theorem \ref{upper_bound_theorem} states that 
$\ind\left( \Mat_N\otimes \mathfrak A\right)$ cannot exceed the power of the
factor $\lambda+\mu$ in characteristic polynomial of $\Mat_N\otimes \mathfrak A$.
This polynomial can be computed using theorem \ref{ext_cayley}:
\begin{eqnarray}
\lefteqn{\chi_{\Mat_N\otimes \mathfrak A}=\det\left((\lambda A+\mu A^T)^N \det(F) \prod_{i\neq j}
{(\lambda\alpha_i A+\mu \alpha_j A^T)}   \right)=}\cr
& &\qquad=\left(\det(\lambda A+\mu A^T)\right)^N \det(F)^{\dim \mathfrak A}
\prod_{i\neq j}{\det(\lambda\alpha_i A+\mu \alpha_j A^T)}=\cr
& &\qquad=\left(\chi_{\mathfrak A}\right)^N\det(F)^{\dim \mathfrak A}
\prod_{i\neq j}{\det(\lambda\alpha_i A+\mu \alpha_j A^T)} \nonumber
\end{eqnarray}
Since we know that $\ind \mathfrak A$ is equal to the power of the factor
$\lambda+\mu$ the term $\left(\chi_{\mathfrak A}\right)^N$ is divisible by
exactly $\lambda+\mu$ to the power $N\ind \mathfrak A$. The other factors do 
not contribute since $\det(F)$ does not contain $\lambda$ or $\mu$ and $\alpha_i$
are independent.
\end{proof}

\section{Symmetric form on $\stab_F$}
\begin{definition}
Let $\mathfrak A$ be an associative algebra. For an element $F\in {\mathfrak A}^*$
 we define
$$
\Stab_F=\left\{a\in {\mathfrak A}:\forall x\in {\mathfrak A}\imply F(ax)=F(xa)\right\}
$$
\end{definition}

\begin{proposition}
\begin{enumerate}
\item $\Stab_F$ is a subalgebra in $\mathfrak A$.
\item If $\mathfrak A$ possesses a unity then $\Stab_F$ possesses a unity.
\end{enumerate}
\end{proposition}

\begin{definition} We define the form $Q_F$ on $\stab_F$ by
the following formula:
$$
Q_F(a,b)=F(ab)
$$
\end{definition}

\begin{theorem}\buzz{Properties of the form $Q_F$} The form $Q_F$ possesses the following
properties:
\begin{enumerate}
\item $Q_F$ is symmetric.
\item $Q_F$ is $\ad{}$ invariant. 
\end{enumerate}
\end{theorem}
\begin{proof} The proof is not difficult:
\begin{enumerate}
\item 
$$
Q_F(a,b)=F(ab)=F(ba)=Q_F(b,a)
$$
\item 
\begin{eqnarray}
\lefteqn{Q_F(\ad{c}a,b)=F((\ad{c}a)b)=F((ca-ac)b)=}\cr
&&=F(cab-acb)=F(abc-acb)=\cr
&&=-F(a(\ad{c}b))=-Q_F(a,\ad{c}b)\nonumber
\end{eqnarray}
\end{enumerate}
\end{proof}

\begin{theorem}
The index of associative algebra $\mathfrak A$ is equal to the maximal power $N$
when $\left(\lambda+\mu\right)^N$ divides $\chi_{\mathfrak A}$ if and only if
the quadratic form $Q_F$ is non-degenerate for a generic $F\in \mathfrak A^*$.
\end{theorem}

{\bf Corollary} An easy source of examples of associative algebras $\mathfrak A$
that satisfy this condition are algebras with unity with index
equal to $1$. Indeed, the dimension of stabilizer of such algebra is $1$ (which is
thus generated by unity) and it is not difficult to find generic $F$ that does
not vanish on the stabilizer (which is the same for all generic $F$).

In particular for these algebras $\ind \left(\Mat_n \otimes \mathfrak A\right)=n$.
\begin{proof}
We will show that non-degeneracy of $Q_F$ for a particular $F$ implies that the
maximal power of $\lambda+\mu$ that divides $\chi(\cdot,\cdot,F)$ (evaluated on $F$)
is equal to $\dim \Stab_F$. This will imply the statement of the theorem in case
$F$ is generic.

Let $A$ be the multiplication table of $\mathfrak A$ evaluated in $F$. Let us choose a basis 
$\left\{e_i\right\}$ in $\Stab_F$ and complement it with vectors $\left\{w_j\right\}$
to the full basis in $\mathfrak A$ in such a way that 
$(A-A^T)$ restricted to $\left\{w_j\right\}$ is a non-degenerate skew-symmetric form.

Then 
$$
\chi(\lambda,\mu,F)=\det(\lambda A+\mu A^T)=\det((\lambda+\mu)A-\mu(A-A^T))
$$

The matrix $A-A^T$ is a block matrix with respect to partition of basis in 
$\left\{e_i\right\}$ and $\left\{w_j\right\}$. The only non-zero block is the
block $\left\{w_j\right\}\times\left\{w_j\right\}$. Thus the only way for the
term 
$$
\left(\lambda+\mu\right)^{\dim \Stab_F}\mu^{\dim \mathfrak A-\dim \Stab_F}
$$
to appear in $\det((\lambda+\mu)A-\mu(A-A^T))$ is when $A$ restricted to 
$\left\{e_i\right\}$ is non-degenerate. But $\left.A\right|_{\left\{e_i\right\}}=Q_F$.
\end{proof}

\references

\end{document}